\documentclass{owrart}
\usepackage{amscd,amssymb}
\usepackage[colorlinks,plainpages,urlcolor=blue]{hyperref}
\usepackage[all]{xy}

\newcommand{\Z}{\mathbb{Z}}
\newcommand{\Q}{\mathbb{Q}}
\newcommand{\R}{\mathbb{R}}
\newcommand{\C}{\mathbb{C}}

\newcommand{\RP}{\mathbb{RP}}
\newcommand{\CP}{\mathbb{CP}}

\newcommand{\B}{\mathcal{B}}
\newcommand{\T}{\mathcal{T}}

\newcommand{\abs}[1]{\left| #1 \right|}

\theoremstyle{plain}

\theoremstyle{definition}

\begin{document}

\begin{talk}{Alexander I. Suciu}
{The rational homology of real toric manifolds}{Suciu, Alex}

\subsection*{Toric manifolds}

In a seminal paper \cite{DJ91} that appeared some twenty 
years ago, Michael Davis and Tadeusz Janusz\-kiewicz 
introduced a topological version of smooth toric varieties, 
and showed that many properties previously discovered 
by means of algebro-geometric techniques are, in fact, 
topological in nature.  

Let $P$ be an $n$-dimensional simple polytope with facets 
$F_1,\dots, F_m$,  and let $\chi$ be an integral $n \times m$ 
matrix such that, for each vertex $v=F_{i_1}\cap \cdots \cap F_{i_n}$, 
the minor of columns $i_1,\dots, i_n$ has determinant 
$\pm 1$.  To such data, there is associated a $2n$-dimensional 
toric manifold, $M_P(\chi)=T^n\times P/\sim$, 
where $(t,p)\sim (u,q)$ if $p=q$, and $tu^{-1}$ belongs to 
the image under $\chi \colon T^m\to T^n$ of the coordinate 
subtorus corresponding to the smallest face of $P$ containing $q$ 
in its interior.  

Here is an alternate description, using the moment-angle complex 
construction (see for instance \cite{DS07} and references therein).  
Given a simplicial complex $K$ on vertex set $[n]=\{1,\dots, n\}$, 
and a pair of spaces $(X,A)$, let $\mathcal{Z}_K(X,A)$ be 
the subspace of the  cartesian product $X^{\times n}$, defined 
as the union $\bigcup_{\sigma\in K}  (X,A)^{\sigma}$, 
where $(X,A)^{\sigma}$ is the set of points for which the 
$i$-th coordinate belongs to $A$, whenever $i\notin \sigma$.  
It turns out that the quasi-toric manifold $M_P(\chi)$ is obtained 
from the moment angle manifold  $\mathcal{Z}_K(D^2,S^1)$, 
where $K$ is the dual to $\partial P$, 
by taking the quotient by the 
relevant free action of the torus $T^{m-n}=\ker(\chi)$.  

\subsection*{Real toric manifolds}

An analogous theory works for real quasi-toric manifolds, 
also known as small covers. Given a homomorphism
$\chi \colon \Z_2^m \to \Z_2^n$ satisfying a minors 
condition as above, the resulting $n$-dimensional 
manifold, $N_P(\chi)$, is the quotient of the real moment 
angle manifold $\mathcal{Z}_K(D^1,S^0)$ by a free 
action of the group $\Z_2^{m-n}=\ker(\chi)$. The manifold 
$N_P(\chi)$ comes equipped with an action of $\Z_2^n$;  
the associated Borel construction is homotopy equivalent 
to $\mathcal{Z}_K(\RP^\infty,*)$. 

If $X$ is a smooth, projective toric variety, 
then $X(\C)=M_P(\chi)$, for some simple polytope 
$P$ and characteristic matrix $\chi$, 
and $X(\R)=N_P(\chi \bmod 2\Z)$. Not all toric 
manifolds arise in this manner. For instance, 
$M=\CP^2 \sharp \CP^2$  is a toric 
manifold over the square, but it does not admit 
any (almost) complex structure; thus, $M\not\cong X(\C)$.

The same goes for real toric manifolds. For instance, 
take $P$ to be the dodecahedron, and use one of 
the characteristic matrices $\chi$ listed in \cite{GS03}. 
Then, by a theorem of Andreev \cite{An70},  the 
small cover $N_P(\chi)$ is a hyperbolic $3$-manifold; 
thus, by a theorem of Delaunay \cite{De05}, 
$N_P(\chi)\not\cong X(\R)$. 

\subsection*{The Betti numbers of real toric manifolds}

In \cite{DJ91}, Davis and Januszkiewicz showed that the 
sequence of mod~$2$ Betti numbers of $N_P(\chi)$ 
coincides with the $h$-vector of $P$. 
In joint work with Alvise Trevisan \cite{ST12}, 
we compute the rational cohomology groups 
(together with their cup-product structure) for real, 
quasi-toric manifolds.  It turns out that the rational 
Betti numbers are much more subtle, depending 
also on the characteristic matrix $\chi$. 

More precisely, for each subset $S\subseteq [n]$, let 
$\chi_S = \sum_{i \in S} \chi_i$, where $\chi_i$ is the 
$i$-th row of $\chi$, and let $K_{\chi,S}$ be the induced 
subcomplex of $K$ on the set of vertices $j \in [m]$ 
for which the $j$-th entry of $\chi_S$ is non-zero. 
Then 
\begin{equation}
\label{eq:betti npchi} \tag{*}
\dim H_q(N_P(\chi),\Q)= \sum_{S\subseteq [n]}
\dim \widetilde{H}_{q-1}(K_{\chi, S},\Q).
\end{equation}

The proof of formula \eqref{eq:betti npchi}, given in \cite{ST12}, 
relies on two fibrations relating the real toric manifold $N_P(\chi)$ 
to some of the aforementioned moment-angle complexes,
\[
\xymatrixrowsep{14pt}
\xymatrixcolsep{10pt}
\xymatrix{
& \Z_2^{m-n} \ar[d] \\
& \mathcal{Z}_K(D^1,S^0) \ar[d] \\
\Z_2^{n} \ar[r]& N_P(\chi)  \ar[r]& \mathcal{Z}_K(\RP^{\infty},*)\, .}
\]
The proof entails a detailed analysis of homology in rank 
$1$ local systems on the space $\mathcal{Z}_K(\RP^\infty,*)$, 
exploiting at some point the stable splitting of moment-angle 
complexes due to Bahri, Bendersky, Cohen, and Gitler \cite{BBCG10}. 
Some of the details of the proof appear in Trevisan's 
Ph.D. thesis \cite{Tr12}. 

As an easy application of  formula \eqref{eq:betti npchi}, one 
can readily recover a result of Nakayama and Nishimura \cite{NN05}:  
A real, $n$-dimensional toric manifold $N_P(\chi)$ is orientable 
if and only if there is a subset $S\subseteq [n]$ such that 
$K_{\chi,S}=K$.

\subsection*{The Hessenberg varieties}

A classical construction associates to each 
Weyl group $W$ a smooth, complex projective toric variety 
$\T_W$, whose fan corresponds to the reflecting hyperplanes 
of $W$ and its weight lattice. 

In the case when $W$ is the symmetric group $S_n$, 
the manifold $\T_n=\T_{S_n}$ is the well-known Hessenberg 
variety, see \cite{DPS92}.  Moreover, $\T_n$ is isomorphic to the 
De~Concini--Procesi wonderful model $\overline{Y_\mathcal{G}}$, 
where $\mathcal{G}$ is the maximal building set for the Boolean 
arrangement in $\CP^{n-1}$.  Thus, $\T_n$ can be obtained 
by iterated blow-ups:  first blow up $\CP^{n-1}$ at the $n$ 
coordinate points, then blow up along the proper transforms 
of the $\binom{n}{2}$ coordinate lines, etc. 

The real locus, $\T_n(\R)$, is 
a smooth, real toric variety of dimension $n-1$; 
its rational cohomology was recently computed by 
Henderson \cite{He12}, who showed that 
\[
\dim H_i(\T_n(\R),\Q) = A_{2i} \binom{n}{2i}, 
\]
where $A_{2i}$ is the Euler secant number, defined as the 
coefficient of $x^{2i}/(2i)!$ in the Maclaurin expansion of $\sec(x)$.
As announced in \cite{Su12}, we can recover this computation, 
using formula \eqref{eq:betti npchi}. 

To start with, note that the $(n-1)$-dimensional polytope 
associated to $\T_n(\R)$ is the permutahedron $P_n$. 
Its vertices are obtained by permuting the coordinates 
of the vector $(1,\dots, n)\in \R^n$, while its facets  
are indexed by the non-empty, proper subsets $Q\subset [n]$. 
The characteristic matrix $\chi = (\chi^Q)$ for $\T_n(\R)$ 
can be described as follows:  $\chi^i$ is the $i$-th 
standard basis vector of $\R^{n-1}$ for $1\le i<n$, while 
$\chi^n=\sum_{i<n} \chi^i$ and $\chi^Q = \sum_{i \in Q} \chi^i$.

The simplicial complex $K_n$ dual to $\partial P_n$ 
is the barycentric subdivision of the boundary of 
the $(n-1)$-simplex.  Given a subset $S\subset [n-1]$, 
the induced subcomplex $(K_n)_{\chi,S}$ depends only 
on the cardinality $r=\abs{S}$; denote any one of these 
$\binom{n-1}{r}$ subcomplexes by $K_{n,r}$.  It turns 
out that $K_{n,r}$ is the order complex associated to 
a rank-selected poset of a certain subposet 
of the Boolean lattice $B_n$. A result of 
Bj\"{o}rner and Wachs \cite{BW} insures that such 
simplicial complexes are Cohen--Macaulay, and 
thus have the homotopy type of a wedge of spheres 
(of a fixed dimension); in fact, 
$K_{n,2r-1}\simeq K_{n,2r} \simeq \bigvee^{A_{2r}} S^{r-1}$. 
Hence,
\begin{align*}
\dim H_i( \T_n(\R),\Q ) &= 
\sum_{S \subseteq [n-1] } \dim \widetilde{H}_{i-1}( (K_n)_{\chi,S},\Q )\\
 &= \sum_{r=1}^{n-1}  \binom{n-1}{r}  \dim \widetilde{H}_{i-1}( K_{n,r},\Q)\\
&= \left( \binom{n-1}{2i-1} + \binom{n-1}{2i}\right) A_{2i} 
 =  \binom{n}{2i} A_{2i}.
\end{align*}

Recently, Choi and Park \cite{CP12} have extended this 
computation to a much wider class of real toric manifolds. 
Given a finite simple graph $\Gamma$, let $\B(\Gamma)$ 
be the building set obtained from the connected induced 
subgraphs of $\Gamma$, and let $P_{\B(\Gamma)}$ be 
the corresponding graph associahedron.   Using formula 
\eqref{eq:betti npchi}, these authors compute the Betti 
numbers of the smooth, real toric variety $X_{\Gamma}(\R)$ 
defined by $P_{\B(\Gamma)}$.  When $\Gamma=K_n$ is 
a complete graph, $X_{K_n}=\T_n$, and one recovers 
the above calculation.
 
\subsection*{The formality question}


 A finite-type CW-complex $X$ is said to be {\em formal}\/ 
 if its Sullivan minimal model is quasi-isomorphic to 
the rational cohomology ring of $X$, endowed 
with the $0$ differential. Under a nilpotency 
assumption, this means that $H^*(X,\Q)$ 
determines the rational homotopy type of $X$.

As shown by Notbohm and Ray \cite{NR05}, 
if $X$ is formal, then $\mathcal{Z}_K(X,*)$ is formal; 
in particular, $\mathcal{Z}_K(S^1,*)$ and  
$\mathcal{Z}_K(\CP^{\infty},*)$ are always formal. 
More generally, as shown by F\'elix and Tanr\'e \cite{FT08}, 
if both $X$ and $A$ are formal, and the inclusion 
$A \hookrightarrow X$ induces a surjection 
in rational cohomology, then 
$\mathcal{Z}_K(X,A)$ is formal. 

On the other hand, as sketched in \cite{Ba03}, and proved 
with full details in \cite{DS07}, the spaces 
$\mathcal{Z}_K(D^2,S^1)$ can have 
non-trivial triple Massey products, and thus are 
not always formal.   In fact, as shown in  \cite{DS07}, 
there exist polytopes $P$ and dual triangulations 
$K=K_{\partial P}$ for which the moment-angle 
manifold $\mathcal{Z}_K(D^2,S^1)$ is not formal. 
Using these results, as well as a construction from 
\cite{BBCG}, we can exhibit real moment-angle 
manifolds $\mathcal{Z}_L(D^1,S^0)$ that are not formal.

In view of this discussion, the following natural question 
arises:  are toric manifolds formal?  Of course, smooth 
(complex) toric varieties are formal, by a classical result 
of Deligne, Griffith, Morgan, and Sullivan.  More generally,  
Panov and Ray showed in \cite{PR08} that all  
toric manifolds are formal.  So we are left with the 
question whether real toric manifolds are always formal.

\subsection*{Acknowledgement} Research partially supported 
by NSA grant H98230-09-1-0021 and NSF grant DMS--1010298.

\newcommand{\arxiv}[1]
{\texttt{\href{http://arxiv.org/abs/#1}{arXiv:#1}}}
\newcommand{\doi}[1]
{\texttt{\href{http://dx.doi.org/#1}{doi:#1}}}

\end{talk}

\end{document}